\input amstex
\documentstyle{amsppt}
\magnification=\magstep1
\hsize=5in
\vsize=7.3in
\TagsOnRight
\topmatter
\title On the signature of certain intersection forms 
\endtitle
\author Feng Xu \endauthor
\abstract We prove a conjecture of Zuber on the signature of intersection
froms associated with affine algebras of type A.

\endabstract
 
\thanks
This work is partially
supported by NSF grant DMS-9500882.\endthanks
\endtopmatter
\heading \S1.  Introduction \endheading
Let $N\geq 2$ be a positive integer and 
$\lambda'= (\lambda_1',...,\lambda_{N-1}'),
0< \lambda_i' < 1, i=1,2,...,N-1$
with $\sum_{0<i<N} \lambda_i' <1$. 
Define
$$
p_i(\lambda') := \lambda_i' + ...+ \lambda_{N-1}' - \frac{1}{N}
\sum_{0<j<N} j\lambda_j', i=1,2,...,N-1
$$ 
and $p_N(\lambda') :=  - \frac{1}{N}\sum_{0<j<N} j\lambda_j'$.
Define 
$$
q_i(\lambda') := -Np_i(\lambda') + \frac{N+1}{2} - i
$$ and
$$
g_i(\lambda') := (-1)^i \prod_{r=1}^{N} 
2cos(\pi (p_r(\lambda') -  \frac{i}{N}))  
$$, $i=1,2,...,N$. \par
If $S$ is a finite sequence of real numbers, we define $b_+(S)$(resp.  $b_-(S),
b_0(S)$) to be the number of positive (resp. negative, zero) elements in 
$S$. Let $a(S):= b_+(S) - b_-(S)$ and denote by $Q_{\lambda'},G_{\lambda'}$
the following two sets:
$$
Q_{\lambda'}:= \{ cos(\pi q_1(\lambda')),..., cos(\pi q_N(\lambda')) \},
$$
$$
G_{\lambda'}:= \{ g_1(\lambda'),...,  g_N(\lambda') \}.
$$
Notice that since $ cos(\pi q_i(\lambda')) >0$ iff $ q_i(\lambda')\in 
]2p-\frac{1}{2}, 2p+\frac{1}{2}[$ for some integer $p$, $b_+(Q_{\lambda'})$
is much easier to calculate than  $b_+(G_{\lambda'})$ and the same is true
for $b_-$'s. The main theorem in this paper is the following:
\proclaim{Theorem 1}
Let $\lambda' = (\lambda'_1,...,\lambda'_{N-1})$ be as above. Then:
$$
b_+(Q_{\lambda'}) = b_+(G_{\lambda'}), 
b_-(Q_{\lambda'}) = b_-(G_{\lambda'}),
b_0(Q_{\lambda'}) = b_0(G_{\lambda'})
$$.
\endproclaim
This theorem implies Zuber's conjecture about
the signature of intersection forms associated
with affine algebras of type A (cf.[Z1]) which is the motivation
of this paper. Note that $b_0(Q_{\lambda'}) = b_0(G_{\lambda'})$
is already noticed in a slightly different form in [Z1]. \par
Zuber's conjecture appeared as Conjecture 2.5 of [Z1]. It is based on the
mysterious connections between integrable models with two supersymmetries
(N=2) in two dimensions (cf.[CV]) and the class of graphs constructed in
[Z1] (see also [Z2]). In the special case when the graphs are regular
(cf.\S2.3), the conjecture can be proved (cf. Page 14 of [GZV]) 
by combining the results
of [GZV] and [S]. In fact, in [GZV], a connection between regular graphs 
and singularity theory is established, and combine with [S] which is
based on mixed Hodge structures, gives a rather indirect proof of
Zuber's conjecture in  the special case when the graphs are regular.
However, for other graphs in [Z1] (see also [X]), the connection with 
singularity theory, if any, is not clear at all.  We also failed to prove 
Zuber's conjecture by using the connection with [CV] as mentioned in
[Z1].

We came to the realization that a statement as in theorem 1 may be
true by first observing lemma 1 (cf.\S2.1) which was already noticed
in a slightly different form in [Z1].  We then checked that theorem 1 is
true explicitly in the case when $N=3,4$ and some other cases which
motivated us to give a general proof.

The idea of the proof of theorem 1 is as follows. When $\lambda'$ changes,
$ b_+(G_{\lambda'})$(resp. $ b_-(G_{\lambda'})$) may change only if
some of $g_i(\lambda')$'s become $0$ or change its sign, i.e., 
$g_i(\lambda')$ intersect the hyperplanes on which 
$g_i(\lambda') =0$. By lemma 1 of \S2.1, these  hyperplanes are the same
as the  hyperplanes on which some  $q_j(\lambda')$'s lie in $\Bbb Z+
\frac{1}{2}$.  Consider the domain $D:= \{ \lambda'=(\lambda'_1,... 
\lambda'_{N-1}) | 0<\lambda'_i<1, \sum_{0<i<N} \lambda'_i <1 \}$. $D$ is
separated by the above  hyperplanes  into disjoint open regions. In each
open region, the set of numbers compared in theorem 1 should be 
completely determined. In \S2.1, we determine these numbers in a given
open region and find that they miraculously satisfy theorem 1. In \S2.2,
we show that theorem 1 also holds for any 
$\lambda' \in D$ which is on  the boundary of the open region: this follows
from \S2.1 and lemma 1.   In \S2.3, after
introducing Zuber's conjecture, we show how theorem 1 implies
that the conjecture is true. In \S3, we present our conclusions and questions.
\heading{2. The proof} \endheading
We shall use the notations of \S1. Recall \par
$D=\{ \lambda'=(\lambda'_1,...,
\lambda'_{N-1})| \lambda'_i >0, \sum_{0<i<N}  \lambda'_i <1 \}$. 
For $ \lambda' \in D$, recall
$$
\align
p_i(\lambda') &= \lambda'_i +...+ \lambda'_{N-1} - \frac{1}{N} \sum_{0<j<N}
j \lambda'_j \\
&= \frac{1}{N} [(-\lambda'_1 -2 \lambda'_2-...-(i-1)\lambda'_{i-1}) +
(N-i)\lambda'_i +...+ \lambda'_{N-1}]
\endalign
$$ for $i=1,2,...,N-1$.
We have:
$$
\align
p_i(\lambda') & > \frac{1}{N}  
[(-\lambda'_1 -2 \lambda'_2-...-(i-1)\lambda'_{i-1})] \\
& > -\frac{1}{N}  (i-1)  \sum_{0<j<N}\lambda'_j >  
-\frac{1}{N}  (i-1)
\endalign
$$, and
$$
\align
p_i(\lambda') & < \frac{1}{N} [(N-i)\lambda'_i +...+ \lambda'_{N-1}] \\
& < \frac{1}{N} (N-i) \sum_{0<j<N} \lambda'_j  < \frac{1}{N} (N-i)
\endalign
$$.
Similarly one can show $ \frac{1}{N} (1-N) < p_N(\lambda')<0$.
So we have:
$$
\frac{1}{N} (1-i) < p_i(\lambda') < \frac{1}{N} (N-i)
$$, $i=1,2,...,N$.
If for some $i$, $ q_i(\lambda') = -N p_i(\lambda') + \frac{1}{2}(N+1) -i
= j+\frac{1}{2}$ with $j\in \Bbb Z$, then 
$$
i-\frac{1}{2}N < i+j < i+\frac{1}{2}N -1
$$.
Define $0<r_{\lambda'} (i)<N+1$ to be the unique integer 
such that $\frac{1}{N} (r_{\lambda'} (i) + j +i) \in \Bbb Z$. In fact,
if $ j +i <0, r_{\lambda'} (i)= -(j+i)$, if $0\leq j +i < N,
 r_{\lambda'} (i)= N-(j+i)$, and if $N\leq (j+i)< \frac{3}{2}N,
 r_{\lambda'} (i)= 2N-(j+i)$.  It follows that 
$$
g_{ r_{\lambda'} (i)}(\lambda') =0
$$.
We have the following lemma which, in a slightly different form also 
appeared on Page 17 of [Z1].
\proclaim{lemma 1}
For any $\lambda'\in D$, the map $i\rightarrow  r_{\lambda'} (i)$
defined above from $\{i| q_i(\lambda') \in \Bbb Z +\frac{1}{2} \}$
to  $\{ r|  g_r(\lambda') =0 \}$ is a one to one and onto map.  Moreover
$ q_i(\lambda') \in  \Bbb Z +\frac{1}{2}$ iff 
$g_{ r_{\lambda'} (i)}(\lambda') =0$ and $ r_{\lambda'} (i)$ depends
only on $i$ and $ q_i(\lambda')$.
\endproclaim
\demo{Proof}
By using definitions we have
$$
g_r(\lambda') = (-1)^r \prod_{i=1}^{N} 2 sin (\frac{\pi}{N}(
q_i(\lambda') +i + r - \frac{1}{2}))
$$. If $g_r(\lambda')=0$, then there exists  $i:=i_{\lambda'} (r)$ such that
$$
\frac{1}{N}(q_i(\lambda') +i + r - \frac{1}{2}) \in \Bbb Z
$$. Let $j\in \Bbb Z$ with $q_i(\lambda') +i + r - \frac{1}{2} = j+i+r$,
then $q_i(\lambda') = -N p_i(\lambda') + \frac{1}{2}(N+1) -i \in 
\Bbb Z + \frac{1}{2}$. Using the fact that $ 0<|p_a - p_b|<1$ for
any $1\leq a\neq b \leq N$, it is 
easy to see that such an  $i:=i_{\lambda'} (r)$ is also unique.  
It is then easy to
check that the map $i\rightarrow  r_{\lambda'} (i)$ and
$r\rightarrow  i_{\lambda'} (r)$ are inverse to each other. The rest
of the lemma follows from the definitions of  $r_{\lambda'} (i)$.  
\enddemo
\hfill Q.E.D.
\par
Let $\beta_i :=  \frac{1}{2}N -i + \gamma_i$ 
with $\gamma_i \in \Bbb Z, i=1,2,...,N$.
Let $D_\gamma := \{\lambda' \in D| \beta_i < N p_i(\lambda') < 
\beta_i +1, i=1,2,...,N \}$.   $D_\gamma$ will be called open regions.
It is clear that if $\gamma \neq \gamma'$, then $D_\gamma \cap D_{\gamma'}
=\emptyset$. Notice that $q_i(\lambda') \in  
\Bbb Z + \frac{1}{2}$ iff $\lambda'$ lies on the boundary of some  
$D_\gamma$, and by lemma 1,  $g_{r_{\lambda'}(i)}(\lambda') =0$ 
iff $\lambda'$ lies on the boundary of some
$D_\gamma$.
Suppose $D_\gamma \neq \emptyset$ and $\lambda' \in D_\gamma$. Then
we have: \par
(1) If $i<j$, then $\beta_i \geq \beta_j$ which follows from the fact
that  $i<j$, then $p_i(\lambda') > p_j(\lambda')$ and 
$\beta_i - \beta_j \in \Bbb Z$; \par
(2) $\beta_1 - \beta_N \leq N$ which follows from  
$p_1(\lambda') - p_N(\lambda') <1$. \par
By (1) we can assume that
$$
\beta_1=...=\beta_{i_1}>\beta_{i_1+1}=...=\beta_{i_2}>...>\beta_{i_{t-1}+1}
=...=\beta_{i_t}
$$
where $1\leq i_1<i_2<...<i_t=N$.
We determine the sign of $g_r(\lambda')$ for a fixed $1\leq r\leq N$.
Since $0<p_1(\lambda') - p_N(\lambda') <1$, there is at most one  
$k+\frac{1}{2}$ with $k\in \Bbb Z$ such that 
$$
p_N(\lambda')-\frac{r}{N} < k+\frac{1}{2} < p_1(\lambda')-\frac{r}{N}
$$. Also notice that if 
$$
\frac{\beta_i -r}{N} < k+\frac{1}{2} <\frac{\beta_i-r+1}{N} 
$$, then we have
$$
\gamma_i -i-Nk <r< \gamma_i -i-Nk+1
$$ which is impossible since $r,\gamma_i$ are integers.  So if there is a
$k\in \Bbb Z$ such that
$$
p_N(\lambda')-\frac{r}{N} < k+\frac{1}{2} < p_1(\lambda')-\frac{r}{N}
$$, then there is a unique integer, denoted by $1\leq f(r)\leq t-1$, such that:
$$
\frac{\beta_{i_{f(r)+1}} +1 -r}{N} \leq k +\frac{1}{2} \leq
\frac{\beta_{i_{f(r)}}  -r}{N}
$$, and the sign of $g_r(\lambda')$ is:
$$
(-1)^r (-1)^{kk'} (-1)^{(k-1)(N-k')} = (-1)^{r+Nk - i_{f(r)}}
$$
where $k':=\{ \beta_j: \beta_j < \beta_{f(r)} \}^\sharp =N-i_{f(r)}
$.   If there is no $k_1\in \Bbb Z$ such that
$$
p_N(\lambda')-\frac{r}{N} < k_1+\frac{1}{2} < p_1(\lambda')-\frac{r}{N}
$$, then there is a $k\in \Bbb Z$ such that:
$$
k-\frac{1}{2} \leq \frac{\beta_N -r}{N}< \frac{\beta_1+1 -r}{N} 
\leq k+\frac{1}{2}
$$, and the sign of  $g_r(\lambda')$ is:
$$
(-1)^{r+kN}
$$. We define $f(r) =t$ in this case.  Let $s:=r+kN$ 
, then the signs of the set
$\{g_r(\lambda')\}$
with $1\leq f(r)\leq t-1$ are given by
$$
(-1)^{s-i_{f(r)}}
$$ with $ \gamma_{i_{f(r)+1}} +1 - i_{f(r)+1} \leq s \leq   
 \gamma_{i_{f(r)}}  - i_{f(r)}$, and the sign of the set
$\{g_r(\lambda')\}$
with $f(r)=t$ is given by
$$ 
(-1)^{s-i_{t}}
$$ with $ \gamma_{i_1} +1 - i_{1}-N \leq s \leq
 \gamma_{i_{t}}  - i_{t}$.
Now we determine the sign of $cos(\pi q_i(\lambda'))$.  Recall
$\beta_i= \frac{N}{2} -i +\gamma_i,  q_i(\lambda') = -N p_i(\lambda')
+ \frac{N+1}{2} -i$ and $\beta_i<N p_i(\lambda')<\beta_i+1$, we have
$$
-\gamma_i - \frac{1}{2} <  q_i(\lambda') < -\gamma_i + \frac{1}{2}
$$. So $cos(\pi q_i(\lambda')) >0$(resp.  $cos(\pi q_i(\lambda')) <0$)
iff $\gamma_i \in 2\Bbb Z$ (resp.  $\gamma_i \in 2\Bbb Z+1$).
Recall from the introduction we have that
for a finite sequence $S$ of real numbers
$a(S) = b_+(S) -  b_-(S)$.  To save some writing for any integer $x$
we define $\{x\} := \frac{1-(-1)^x}{2}$.  
Then the $a$ of the following sequence
$\{ cos (\pi  q_i(\lambda')), i_{u-1} +1 \leq i \leq i_u \}$
is
$$
(-1)^{ \gamma_{i_{u}}} \{  i_{u}- i_{u-1} \}
$$
and  the $a$ of the following
sequence
$\{ (-1)^{s-i_u},  
\gamma_{i_{u+1}} +1 - i_{u+1} \leq s \leq
 \gamma_{i_{u}}  - i_{u} \}
$ is
$$
(-1)^{ \gamma_{i_{u}}} \{ \gamma_{i_{u}}-\gamma_{i_{u+1}} + i_{u+1}- i_{u} \}
$$, where we define $i_{u-1}=0$ if $u=1$, and $\gamma_{i_{u+1}} -i_{u+1}
=\gamma_{i_{1}} - {i_{1}} -N$ if $u=t$.
It follows that $a(G_{\lambda'}) - a(Q_{\lambda'})$ is given by
$$
\align
\sum_{u=1}^{t} &
(-1)^{ \gamma_{i_{u}}} (\{ \gamma_{i_{u}}-\gamma_{i_{u+1}} + i_{u+1}- i_{u} \}
-  \{  i_{u}- i_{u-1} \}) \\
=& (-1)^{ \gamma_{i_{1}}} (\{ \gamma_{i_{1}}-\gamma_{i_{2}} + i_{2}- 
i_{1} \}
-  \{  i_{1} \})+ \\
& (-1)^{ \gamma_{i_{2}}} (\{ \gamma_{i_{2}}-\gamma_{i_{3}} + i_{3}- i_{2} \}
-  \{  i_{2}- i_{1} \}) + \\
&...+   (-1)^{ \gamma_{i_{t-1}}} (\{ \gamma_{i_{t-1}}-\gamma_{i_{t}} + i_{t}- 
i_{t-1} \}
-  \{  i_{t-1}- i_{t-2} \}) + \\
& (-1)^{ \gamma_{i_{t}}} (\{ \gamma_{i_{t}}-\gamma_{i_{1}} + i_{1} \}
-  \{  i_{t}- i_{t-1} \})
\endalign
$$. By using
$$
\{ \gamma_{i_{u}}-\gamma_{i_{u+1}} + i_{u+1}- i_{u} \}
= \{ \gamma_{i_{u}}-\gamma_{i_{u+1}} \} + 
(-1)^{\gamma_{i_{u}}-\gamma_{i_{u+1}}} \{  i_{u+1}- i_{u} \}
$$ which follows easily from the definition of $\{.\}$, we see that 
$\pm \{  i_{u+1}- i_{u} \}$ terms cancelled each other in the above 
summation and the remaining terms are:
$$
 (-1)^{ \gamma_{i_{1}}} \{ \gamma_{i_{1}}-\gamma_{i_{2}} \} +
 (-1)^{ \gamma_{i_{2}}} \{ \gamma_{i_{2}}-\gamma_{i_{3}} \}
+...+  (-1)^{ \gamma_{i_{t}}} \{ \gamma_{i_{t}}-\gamma_{i_{1}} \}
$$ which is also $0$ since $ \{x \} = \frac{1-(-1)^x}{2}$.
So we have shown that 
$$
a(G_{\lambda'}) - a(Q_{\lambda'})=0
$$, i.e., 
$$
b_+(G_{\lambda'}) - b_-(G_{\lambda'}) =
b_+(Q_{\lambda'}) - b_-(Q_{\lambda'})
$$.  Since
$$
b_+(G_{\lambda'}) + b_-(G_{\lambda'}) = N =  
b_+(Q_{\lambda'}) + b_-(Q_{\lambda'})
$$, it follows that theorem 1 is true for $\lambda'\in D_\gamma$.
\subheading {2.2. The boundary case}
Assume $\lambda'\in D$ and $\lambda'$ is on the boundary of some
$ D_\gamma$. Assume $\{q_i(\lambda') | q_i(\lambda') \in \Bbb Z+ \frac{1}{2}
\} = \{ q_{k_1}(\lambda'),...,q_{k_s}(\lambda'), s\geq 1 \}$.  Let
$k_i\rightarrow r_{\lambda'}(k_i)$ be as in lemma 1.  We can choose
a small neighborhood $W$ of $\lambda'$ such that for any $\mu\in W$ and
$l\neq k_i, i=1,...,s$ (resp. $m\neq r( k_i), i=1,...,s$), 
$cos (\pi q_l(\mu))$ (resp. $g_m(\mu)$) has the same sign as $cos (\pi
 q_l(\lambda'))$
(resp. $g_m(\lambda')$) since  $cos (\pi
 q_l(\lambda'))$ (resp. $g_m(\lambda')$) is not zero.   
Let $\mu_1 \in D_\gamma \cap W$. We compare $b_+(G_{\lambda'})$(resp. 
 $b_+(Q_{\lambda'})$) with  $b_+(G_{\mu_1})$ (resp.  $b_+(Q_{\mu_1})$).
Since $b_+(G_{\mu_1}) = b_+(Q_{\mu_1})$ by \S2.1, to prove
$b_+(G_{\lambda'}) = b_+(Q_{\lambda'})$ we just have to show that 
if $cos (\pi q_{k_i}(\mu_1)) >0$ for some $k_i, 1\leq i\leq s$, then
$g_{r_{\lambda'}(k_i)}(\mu_1) >0$ and vice versa.
Let us consider a small line segment with end points
$\mu_1, \mu_2$  which passes from $D_\gamma$ to its neighbor $D_{\gamma'}$, 
intersects the hyperplane $q_{k_i}(\mu) = q_{k_i}(\lambda')$ at $\mu_0$,
and does not intersect any other hyperplanes.  Then we have 
$cos(\pi q_{k_i}(\mu_0)) = 0$, so by lemma 1, $g_{r_{\mu_0}(k_i)}(\mu_0) = 0$
. Again by lemma 1, $r_{\mu_0}(k_i)$ depends only on $k_i$ and 
$q_{k_i}(\mu_0) = q_{k_i}(\lambda')$, so  $r_{\mu_0}(k_i)= r_{\lambda'}(k_i)$.
As $\mu$ goes from $\mu_1$ to $\mu_2$ on the above line
segment, 
$cos(\pi q_{k_i}(\mu)), g_{r_{\lambda'}(k_i)}(\mu)$ change their signs
while the signs of all other 
$cos(\pi q_{i}(\mu)), g_{j}(\mu)$'s do not change. By \S2.1, 
$b_+(Q_{\mu_l}) = b_+(Q_{\mu_l}), l=1,2$, it follows that 
if $cos (\pi q_{k_i}(\mu_1)) >0$ for some $k_i, 1\leq i\leq s$, then
$g_{r_{\lambda'}(k_i)}(\mu_1) >0$ and vice versa.
So we have proved that
$b_+(G_{\lambda'}) = b_+(Q_{\lambda'})$, and since 
$b_0(G_{\lambda'}) = b_0(Q_{\lambda'})$ by lemma 1, and both $G_{\lambda'}$ and
 $Q_{\lambda'}$ have $N$ elements, theorem 1 is proved for $\lambda' \in D$
which lies on the boundary of some $D_\gamma$. \par
By \S2.1, \S2.2, theorem 1 is proved.
\subheading {2.3. Zuber's conjecture}
To describe  Zuber's conjecture, we have to introduce some notations from
[Z1] to which the reader is referred for more details.\par
Let $\Lambda_1,...,\Lambda_{N-1}$ be the fundamental weights of $SL(N)$.
Let $k\in \Bbb N$. Recall that the set of integrable weights of the affine
algebra $\widehat{SL(N)}$ at level $k$ is the following subset of the weight
lattice of  $SL(N)$:
$$
P_{++}^{(h)} = \{ \lambda = \lambda_1 \Lambda_1 +...+ \lambda_{N-1}
 \Lambda_{N-1} | \lambda_i\in \Bbb N,  \lambda_1+...+ \lambda_{N-1} <h
\}
$$ where $h=k+N$.  This set admits a $\Bbb {Z}_N$ automorphism generated
by
$$
\sigma: \lambda=( \lambda_1, \lambda_2,..., \lambda_{N-1}) \rightarrow
\sigma( \lambda) = (h-\sum_{j=1}^{N-1} \lambda_j, \lambda_1,...,\lambda_{N-2})
$$.  We then introduce the weights $e_i$ of the standard $N$-dimensional
representation of $SL(N)$
$$
e_1=\Lambda_1,e_i=\Lambda_i - \Lambda_{i-1}, i=2,...,N-1, e_N=-\lambda_{N-1}
$$ endowed with the scalar product $(e_i,e_j)= \delta_{ij} - \frac{1}{N}$.
We shall be concerned with type II class of graphs introduced in
section 1 of [Z1]. These graphs generalize the classical $A,D,E$ Dynkin
diagrams which may be regarded as related to the $SL(2)$ algebra.
The axioms on these graphs are given in \S1.2 of [Z1] as follows: \par
(1) A set $\nu$ of $|\nu|=n$ vertices is given. These vertices are denoted
by Latin letters $a,b,....$ There exists an involution $a\rightarrow \bar{a}$
and the set $\nu$ admits a $\Bbb{Z}_n$ grading denoted by $\tau(a)$ such that
$\tau(\bar{a}) = -\tau(a)$mod$N$;\par
(2) A set of $N-1$ commuting $n\times n$ matrices $G_p,p=1,2,...,N-1$ is 
given.  Their matrix elements are assumed to be non-negative integers, so
they may be regarded as adjacency matrices of $N-1$ graphs {\it g}$_p$.
{\it g}$_1$ is also assumed to be connected;\par
(3) The edges of the graphs  {\it g}$_p$ are compatible with the
grading $\tau$ in the sense that $(G_p)_{ab} =0$ if $\tau(b)\neq \tau(a) +p
$mod$N$;\par
(4) The matrices are transposed of one another
$G_p^{t} = G_{N-p}$ and $(G_p)_{ab} = (G_p)_{\bar{b} \bar{a}}$;\par
(5) As a consequence of (2) and (4), the matrices $G_p$ are commuting normal
matrices and may thus be simultaneously diagonalized in a common orthonormal
basis. This basis, denoted by $\psi^{(\lambda,i)}$, is assumed to be labelled 
by the weights $\lambda$ of $SL(N)$, that are restricted to $P_{++}^{(h)}$,
for some integer $h>N$, in a way that the eigenvalues $\gamma_p^{(\lambda)}$
have the form $ \gamma_p^{(\lambda)} = \chi_p(M(\lambda))$, where 
$ \chi_p$ is the ordinary character for the $p$-th fundamental
representation of the group $SU(N)$, and $M(\lambda)$ denotes the 
diagonal matrix $M(\lambda)=$diag$(\epsilon_j(\lambda))_{j=1,...,N}$.
Here $\epsilon_j(\lambda):= \exp(-\frac{2\pi i}{h}(e_j, \lambda)$, and $i$
in $(\lambda,i)$ is an index integer, $1\leq i\leq m_\lambda$ with 
$m_\lambda$ being the multiplicity of eigenvalue $\gamma_p^{(\lambda)}$.
The set of $(\lambda,i)$'s will be denoted by $Exp$.\par
There exists a special class of solutions known for all $N$ and $h>N$, 
namely the fusion graphs of the affine algebra $\widehat{SL(N)}$ at level
$k=h-N$. The vertices are the integrable weights described above, i.e.,
$\nu = P_{++}^{(h)}$. The matrices $G_p$ are the Verlinde matrices, which
describe the fusion by the $p$-th fundamental representation.  The fusion
rules are given on Page 288 of [K].  Their diagonalization is
known, thanks to the Verlinde formula (cf.  Page 288 of [K]), and the
eigenvalues are the  $\gamma_p^{(\lambda)}$, where $\lambda$ takes all the
values in $= P_{++}^{(h)}$.  We will call these graphs as {\it regular graphs}
in this paper. In the case of $N=2$, these regular graphs reduce to the
$A_{h-1}$ Dynkin diagrams. \par
More solutions are known (cf.[Z1]). In [X] (in 
particular Th.3.10 and (5) of Th.3.8), infinite series of such graphs
are constructed from the maximal conformal inclusions of the form 
$SU(N)\subset G$ with $G$ being a simple and simply connected compact
Lie group.\par
Given graphs of the previous type, let $V$ be a complex vector space
with a basis $\alpha_a$ labelled by the vertices of the set $\nu$. A bilinear
form $g$ is defined by :
$$
g_{ab}= \langle \alpha_a, \alpha_b \rangle = 2\delta_{ab}+ G_{ab}
$$ in terms of the matrix $G_{ab} = \sum_{p=1}^{N-1} (G_p)_{ab}$.  $g$
will be called the intersection form. This is the intersection form 
in the title of this paper. \par
The eigenvalue of the matrix $(g_{ab})$ with eigenvector
$\psi_p^{(\lambda,i)}$ is (cf. (34) of [Z1]):
$$
g^{(\lambda)} = \prod_{i=1}^N (1+ \exp(-\frac{2\pi i}{h}(e_i,\lambda)))
$$.
For  $(\lambda,i) \in Exp$ define real numbers which depend only on
$\lambda$ (cf. (46) of [Z1])by:
$$
q_{\lambda}^{(R)} := \frac{1}{h} \sum_{j=1}^{N-1} j (\lambda_j -1)
+ \frac{(N-h)(N-1)}{2h}
$$.
We can now state Zuber's conjecture on the signature of $g$ (cf. Conjecture
2.5 of [Z1]):
\proclaim{Zuber's Conjecture}
The signature of the bilinear form $g$ for class II graphs is
$(x+, y-, z0)$ where $x$ is the number of 
$q_{\lambda}^{(R)}$ which fall in an interval $]2p-\frac{1}{2}, 
2p+\frac{1}{2}[$ for some $p\in \Bbb Z$ ($p$ may depend on 
$q_{\lambda}^{(R)}$), $y$ is the number of those in an interval
 $]2p'+\frac{1}{2},
2p'+\frac{3}{2}[$ for some $p'\in \Bbb Z$ ($p'$ may depend on
$q_{\lambda}^{(R)}$), and $t=n-r-s$ is the number of those
$q_{\lambda}^{(R)}$ which are half-integers.
\endproclaim
We now prove this conjecture. \par
Let us first notice a simple consequence of the axioms on the graphs.  
It follows from Prop.1.2 of [Z1] that $Exp$ is invariant under the
action of $\sigma$.  In fact, if $\sum_a \psi_a^{(\lambda,i)} a$
is an eigenvector of $G_p$ with eigenvalue $\gamma_p^{(\lambda)}$
, then Prop.1.2 of  [Z1] implies that 
 $\sum_a \psi_a^{(\lambda,i)} \exp(\frac{2\pi i \tau(a)}{N}) a$ is
an eigenvector of $G_p$ with eigenvalue $\gamma_p^{(\sigma(\lambda))}$.
Since $a\rightarrow  \exp(\frac{2\pi i \tau(a)}{N}) a$ is an invertible map,
it follows that  the multiplicity of  
 eigenvalue $\gamma_p^{(\sigma(\lambda))}$ is the same as that of
 eigenvalue $\gamma_p^{(\lambda)}$.  We can therefore define 
$\sigma (\lambda,i) = (\sigma (\lambda),i)$.  It follows that $Exp$ can
be written as a disjoint union of the orbits under the action of
$\sigma$.  To prove Zuber's conjecture, we just have to show 
 it is true on each orbit.  \par
Let $ (\lambda,i) \in Exp$ and
let $d$ be the smallest positive integer such that $\sigma^d (\lambda) = 
\lambda$.  Then $d|N$ and let $N=dd_1$. 
Let $G_\lambda':=\{ g^{(\sigma^i(\mu))}, i=1,2,...,d \}$, 
$Q_\lambda':=\{ cos(\pi q_{\sigma^i(\mu)}^{(R)}), i=1,2,...,d \}$,
we need to show 
$$
b_+(G_\lambda') = b_+(Q_\lambda'),
b_0(G_\lambda') = b_0(Q_\lambda'),
b_-(G_\lambda') = b_-(Q_\lambda').
$$
Note that $b_0(G_\lambda') = b_0(Q_\lambda')$ was already noticed on page
17 of [Z1].
Let $\lambda'= (\frac{\lambda_1}{h},...,\frac{\lambda_{N-1}}{h})$, then
$p_i(\lambda') = \frac{(e_i,\lambda)}{h},i=1,2,...,N$.  
To use theorem 1, we make use of
the following identities which follow from the definitions:
$$
q_{\sigma^{-j}(\lambda)}^{(R)} = q_j(\lambda'), 
$$ and
$$
\align
g^{(\sigma^{j}(\lambda))} &= \prod_{l=1}^{N} (1+ \exp(
+\frac{2\pi i j}{N} \epsilon_l(\lambda))) \\
&= \prod_{l=1}^{N} 2cos(\pi ( p_l(\lambda') - \frac{j}{N})) \times
\exp(\pi ij- \sum_{l=1}^{N} \frac{1}{h} (\lambda,e_l)) \\
&= (-1)^j  \prod_{l=1}^{N} 2cos(\pi ( p_l(\lambda') - \frac{j}{N})) \\
&= g_j(\lambda')
\endalign
$$.
Now it is clear that the size of $G_{\lambda'}$ (resp.  $Q_{\lambda'}$)
is $d_1$ times the size of  $G_{\lambda}'$   (resp.  $Q_{\lambda}'$)
and we have:
$$
b_+(G_{\lambda'}) = d_1 b_+(G_\lambda'),
b_0(G_{\lambda'}) = d_1 b_0(G_\lambda'),
b_-(G_{\lambda'}) = d_1 b_-(G_\lambda')
$$ and
$$
b_+(Q_{\lambda'}) = d_1 b_+(Q_\lambda'),
b_0(Q_{\lambda'}) = d_1 b_0(Q_\lambda'),
b_-(Q_{\lambda'}) = d_1 b_-(Q_\lambda')
$$. By theorem 1, we have proved:
$$
b_+(G_\lambda') = b_+(Q_\lambda'),
b_0(G_\lambda') = b_0(Q_\lambda'),
b_-(G_\lambda') = b_-(Q_\lambda').
$$
Let us summarize the result in the following:
\proclaim{Corollary 1}
Zuber's Conjecture as stated above is true.
\endproclaim
\heading \S 5. Conclusions and questions \endheading
In this paper we proved Zuber's conjecture on the signature of certain
intersection forms by using theorem 1. \par
Our results imply that the infinite series of graphs which are constructed
in [X] by using subfactors associated with conformal inclusions
satisfy Zuber's conjecture.  This lends further support to the idea that these
graphs may be associated with the integrable models in [CV] which
is the basis of  Zuber's conjecture.  Such a relation is not very clear
and should be very interesting.
\heading References \endheading
\roster
\item"{[CV]}"  S.Cecotti and C.Vafa , Nucl.Phys. {\bf 367} (1991) 359-461;
Comm.Math.Phys. {\bf 158} (1993) 569-644.
\item"{[K]}" V.Kac, {\it Infinite dimensional Lie algebras},third edition,
Cambridge University Press.
\item"{[S]}" J.Steenbrink, Composito Math. {\bf 34}(1977) 211-223.
\item"{[GZV]}" S.M. Gusein-Zade and A. Varchenko, Selecta Math. (N.S.) 
{\bf 3}(1997), no.1, 79-97.
\item"{[X]}"  F. Xu, {\it New braided endomorphisms from conformal
inclusions},  \par
Comm. Math. Phys., (1998) (in press).
\item"{[Z1]}"  J.-B. Zuber, {\it Generalized Dynkin diagrams and root 
systems and their folding }, hep-th 9707046, to appear in the
proceedings of the Taniguchi Symposium {\it 
Topological Field theory, Primitive forms and Related Topics,}
Kyoto Dec. 1996, M. Kashiwara, A. Matsuo, K. Saito and I. Satake eds,
Birkha\"{u}ser.
\item"{[Z2]}" J.-B. Zuber, Comm.Math.Phys. {\bf 179} (1996) 265-294.
\endroster
\enddocument